\documentclass{amsart}
\usepackage{amssymb}
\usepackage[dvips]{graphicx}

\newcommand{\gf}{\varphi} \newcommand{\eps}{\varepsilon}
\newcommand{\gl}{\lambda} \newcommand{\ga}{\alpha}
\newcommand{\gb}{\beta}  \newcommand{\gd}{\delta}
\newcommand{\gp}{\psi}   
 
\renewcommand{\gg}{\gamma}  \newcommand{\gs}{\sigma}

 \newcommand{\GP}{\varPsi}
\newcommand{\GF}{\varPhi}   
\newcommand{\GG}{\varGamma} 
\newcommand{\GS}{\varSigma}

\newcommand{\<}{\le} \renewcommand{\>}{\ge}

\newcommand{\cal}{\mathcal}
\renewcommand{\bar}{\overline}
\newcommand{\dbarchi}{\overline{\overline\chi}}
\newcommand{\supp}{\mathop{{\rm supp}}}
\renewcommand{\tilde}{\widetilde}

\newcommand{\bb}{ \mathbb}
\newcommand{\rr}{{\bb R}}
\newcommand{\cc}{{\bb C}}

\newcommand{\z}{{\bb Z}}
\newcommand{\cin}{C^\infty}

\newcommand{\sumi}{\sum^\infty}

\newtheorem{Th}{Theorem}[section] 
\newtheorem{Lem}[Th]{Lemma} 
\newtheorem{Cor}[Th]{Corollary}
\newtheorem{Prop}[Th]{Proposition} 
\theoremstyle{definition}

\theoremstyle{remark}

\numberwithin{equation}{section}
\numberwithin{figure}{section}

\newcommand{\rf}[1]{{\rm(\ref{#1})}}

\newcommand{\sxy}{S''_{xy}}
\newcommand{\ltwor}{L^2(\rr)}
\newcommand{\nt}{\|T\|}
\newcommand{\zp}{\z_+}
\newcommand{\rrp}{\rr_+}

\renewcommand{\)}{)\!)}
\newcommand{\cinxy}{\cin\(x,y\)}
\newcommand{\tup}[1]{\textup{(#1)}}
\newcommand{\D}{\partial}
\renewcommand{\Re}{{\rm Re\,}}
\renewcommand{\Im}{{\rm Im\,}}
\newcommand{\dist}{{\rm dist}}
\begin{document}

\title[Oscillatory integral operators]{Sharp $L^2$ bounds for
oscillatory integral operators with $\cin$ phases}
\author[V.~S.~Rychkov]{Vyacheslav S. Rychkov}
\address{Princeton University, Mathematics Department, Princeton, N.J. 08544, U.S.A.}
\email{rytchkov@math.princeton.edu}
\keywords{Fourier integral operators, Newton polygon.\\
\hphantom{aaa}Supported in part by RFFI grant 99-01-00868.}
\maketitle

\section{Introduction}

Consider the oscillatory integral operator on $\ltwor$ of the form
\begin{equation}
\label{operator}
Tf(x)=\int_{-\infty}^{\infty}e^{i\gl S(x,y)}\chi(x,y)f(y)\,dy 
\end{equation}
with a $\cin$ real phase $S(x,y)$ and a $\cin$ cut-off $\chi(x,y)$
compactly supported in a small neighborhood of the origin in $\rr^2$. We are
interested
in the decay rate of the norm of $T$ on $\ltwor$ for $\gl\to\infty$.

A well-known result of H\"ormander \cite{Hor} says that if the mixed partial derivative 
$\sxy\ne0$ on the support of
$\chi$, then the best possible estimate $\nt\<C{\gl}^{-1/2}$ is true.
The problem of finding the optimal decay rate of $\nt$ for vanishing
$\sxy$ has recently attracted much attention, especially
because of its connection with the smoothing properties of generalized
Radon transforms.

To describe the known results, we need some definitions. Let
$\zp\subset\rrp$
be the sets of non-negative integers and reals. The {\it Newton
polygon of a set} $K\subset \zp^2$ is defined as the convex hull in
$\rrp^2$
of the set $\bigcup_{n\in K}(n+\rrp^2)$. 

Let 
\begin{equation}
\label{taylor}
\sxy(x,y)\sim\sum_{n\in \zp^2}c_n x^{n_1}y^{n_2}
\end{equation}
be the formal Taylor expansion of $\sxy$ at the origin. 
 The Newton
polygon of the set $\{n\in\zp^2|c_n\ne0\}$ is called 
the {\it Newton polygon of the function}  
$\sxy$ and is denoted $\GG(\sxy)$.

Assume that $\GG(\sxy)$ is not empty, which means that the formal
Taylor
series of $\sxy$ is not the zero series. Denote by $t_0$ the parameter
of intersection of the line $n_1=n_2=t$ with the boundary of
$\GG(\sxy)$. The number $\gd=1/(1+t_0)$ is called the {\it 
Newton decay rate} of $S(x,y)$.

This quantity was introduced by Phong and Stein \cite{PS}, who
realized that the decay estimates for $T$ depend on the Newton polygon
of $\sxy$ and proved the bound
\begin{equation}
\label{Phong-Stein}
\nt\<C\gl^{-\gd/2} 
\end{equation}
under the additional assumption that $S(x,y)$ is real-analytic. They
also
showed that this bound is sharp in the sense that 
\[
\nt\>C'\gl^{-\gd/2}
\]
for large $\gl$ if $\chi(0,0)\ne0$, and this part of their result does not depend on
real analyticity. 

On the other hand, the best known estimate in the general 
$\cin$ case was
\begin{equation}
\label{Seeger}
\nt\<C_{\eps}\gl^{-\gd/2+\eps},
\end{equation}
which is implicitly contained in Seeger \cite{See1},\cite{See2}.
It is interesting that the estimates \rf{Phong-Stein} and \rf{Seeger}
have been obtained by quite different arguments.

The purpose of our paper if to show that there no loss
of $\eps$ in the $\cin$ case with one possible exception, when one
loses at most a power of log. 

We say that $\sxy$ is \emph{completely degenerate}, if its formal
Taylor
series \rf{taylor} factorizes in the ring of formal power series 
$\rr[[x,y]]$ as
\[
U(x,y)(y-f(x))^N,
\]
where $N\>2$, the series
$f(x)\in\rr[[x]]$ is of the form $f(x)=cx+\ldots$ with $c\ne0$, and
the series $U(x,y)\in\rr[[x,y]]$ is invertible. 
Note that $\gd=\frac2{N+2}$ for such a phase function. Then we have

\begin{Th} 
\label{main}
There exists a small
neighborhood
of the origin $V$ such that\\
\textup{(a)}
If $\sxy$ is not completely degenerate, and $\supp\chi\subset V$, 
then $\nt\<C\gl^{-\gd/2}$.\\
\textup{(b)} If $\sxy$ is completely degenerate, and $\supp\chi\subset V$,
then 
\begin{equation}
\label{main1}
\nt\<C\gl^{-\frac1{N+2}}(\log\gl)^{\frac{2N}{N+2}}
\end{equation}
\end{Th}
The rest of the paper is devoted to the proof of this theorem, which
is based on a mixture of ideas from the above-mentioned works of Phong and Stein, and
Seeger. Let us describe the scheme of the proof. 

The geometry of the problem is best understood in terms of the
singular variety $\cal Z=\{(x,y)|\sxy=0\}$. Similarly to what Phong
and Stein have done in the real analytic case, we parameterize $\cal
Z$ by means of asymptotic Puiseux series (Section~2). To be more exact, 
assume for the purposes of this discussion that $\sxy$ has the form of a
polynomial in $y$ with $\cin$ coefficients depending on $x$. (In fact,
this case already contains all the essential difficulties.) Then we prove
that for small $x$ there exists a factorization 
\begin{equation}
\label{*fac}
\sxy(x,y)=\prod_{i=1}^n(y-y_i(x)),
\end{equation}
where $y_i(x)$ are continuous $\cc$-valued functions having asymptotic
fractional power series expansions at zero. Now $\cal Z$ splits into
$n$ branches $y=y_i(x)$ (see Fig.~1.1). In this discussion we assume for simplicity
that all these branches are real.

To estimate the norm of $T$, we decompose the operator into pieces 
corresponding to the decomposition of the $(x,y)$ plane into dyadic 
rectangles, so that $T=\sum T_{jk}$, where the kernel of $T_{jk}$ is
supported in $x\approx 2^{-j}$, $y\approx 2^{-k}$.

\begin{tabular}{cc}
\scalebox{0.5}{
\rotatebox{-90}{
\includegraphics{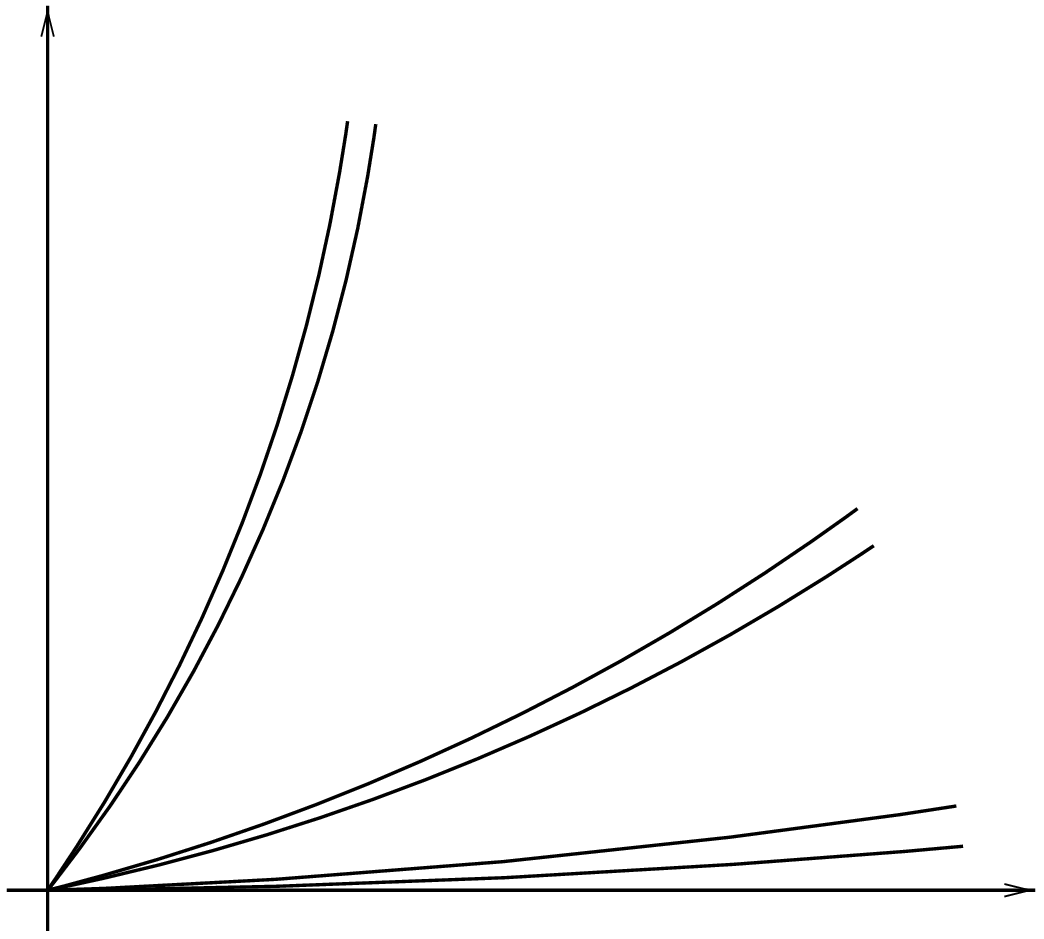}
}}
&
\scalebox{0.5}{
\rotatebox{-90}{
\includegraphics{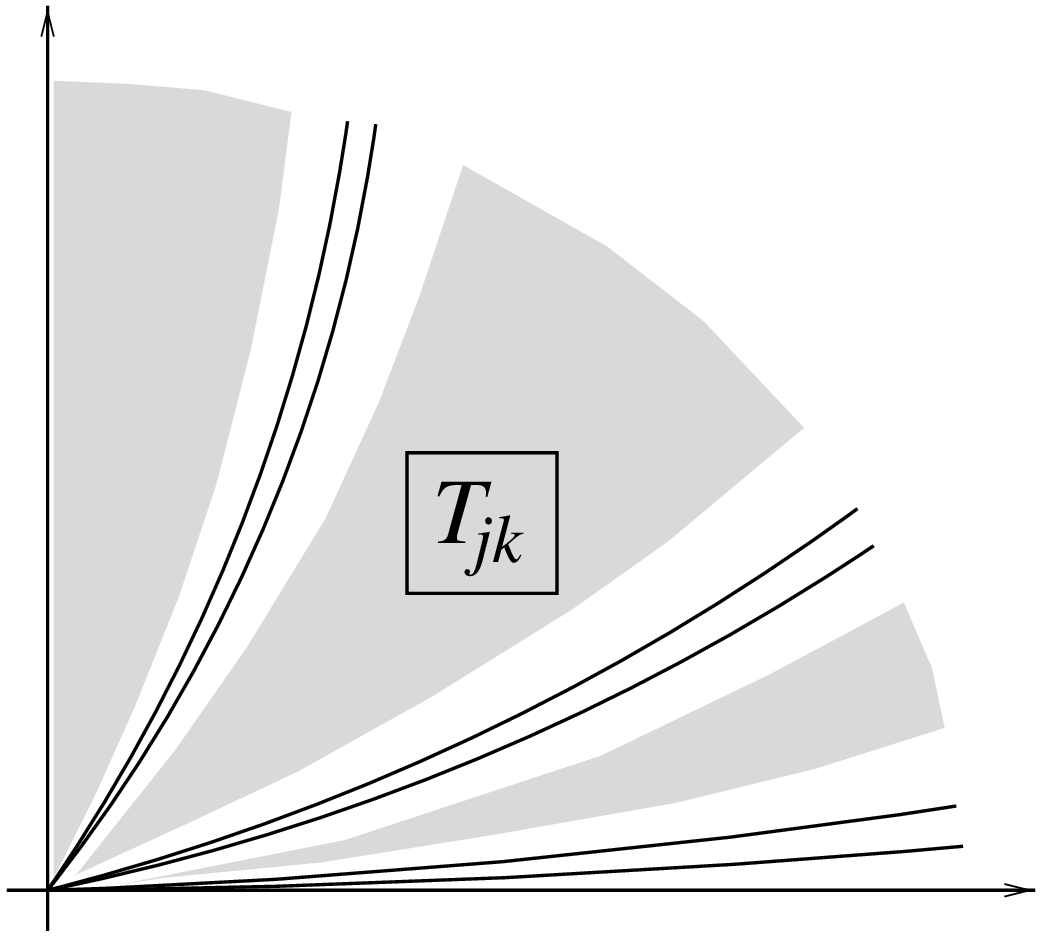}
}}\\
\sc\figurename\ 1.1&
\sc\figurename\ 1.2
\end{tabular}

The way we treat a particular operator $T_{jk}$ depends on its
position with respect to $\cal Z$. The most simple and standard case
is the one of all $T_{jk}$ such that the distance from $\cal Z$ to $\supp T_{jk}$
is relaively large compared to the size of the support. The union of supports of
such $T_{jk}$ comprises the shaded regions in Fig.~1.2.  
For each of these $T_{jk}$ we have from \rf{*fac} a good lower bound
on $\sxy$ on its support, which gives a good estimate on
$\|T_{jk}\|$. Then the required norm estimate for the sum $\sum T_{jk}$ 
over all such $T_{jk}$ is obtained
with the help of a resummation procedure of Phong and Stein, which is based on
almost orthogonality considerations (Section 3).

After that we are left with the $T_{jk}$ supported near $\cal Z$. In
Section 4 we estimate the contribution $T_x=\sum T_{jk}$ of the
operators $T_{jk}$ supported near the branches of $\cal Z$ which are
infinitely tangent to the $x$-axis (the corresponding functions
$y_i(x)$ have the zero asymptotic expansion). Assume that there are
exactly $N$ such branches. We represent $T_x$ as $T_x=\sum T_j$, where
$T_j=\sum_k T_{jk}$ (see Fig.~1.3). Taking into account the fact that the
derivative $\D_y^N\sxy$ does not vanish on the support of $T_j$, it is
easy to obtain by van der Corput's lemma that the kernel of   
the operator $T_j T_j^*$ has the estimate
\[
|K(x,y)|\<C(\gl|x-y|)^{-1/(N+1)}.
\]
From this we can conclude that $\|T_j\|\<C\gl^{-\gd/2}$. Additional
considerations show that the sum $\sum T_j$ is almost orthogonal in
the sense that $\|T_jT_{j'}\|\<C\gl^{-\gd}2^{-\eps|j-j'|}$. Now
the required estimate for $\|T_x\|$ follows by the Cotlar-Stein lemma.

\begin{tabular}{cc}
\scalebox{0.5}{
\rotatebox{-90}{
\includegraphics{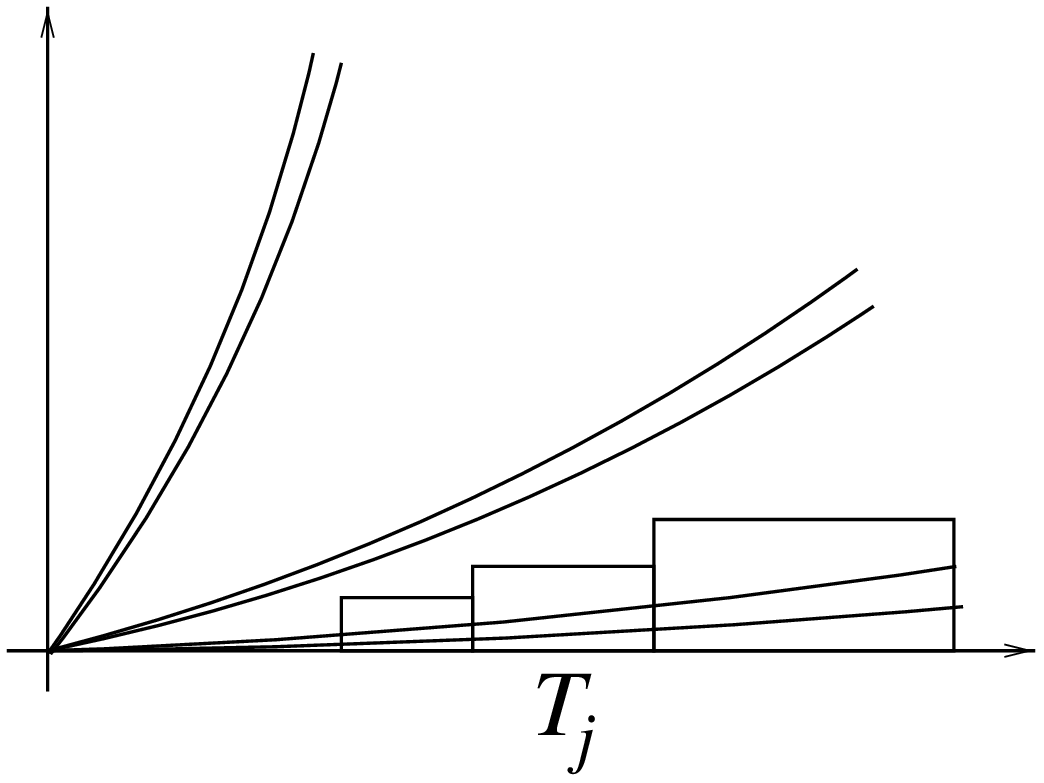}
}}
&
\scalebox{0.5}{
\rotatebox{-90}{
\includegraphics{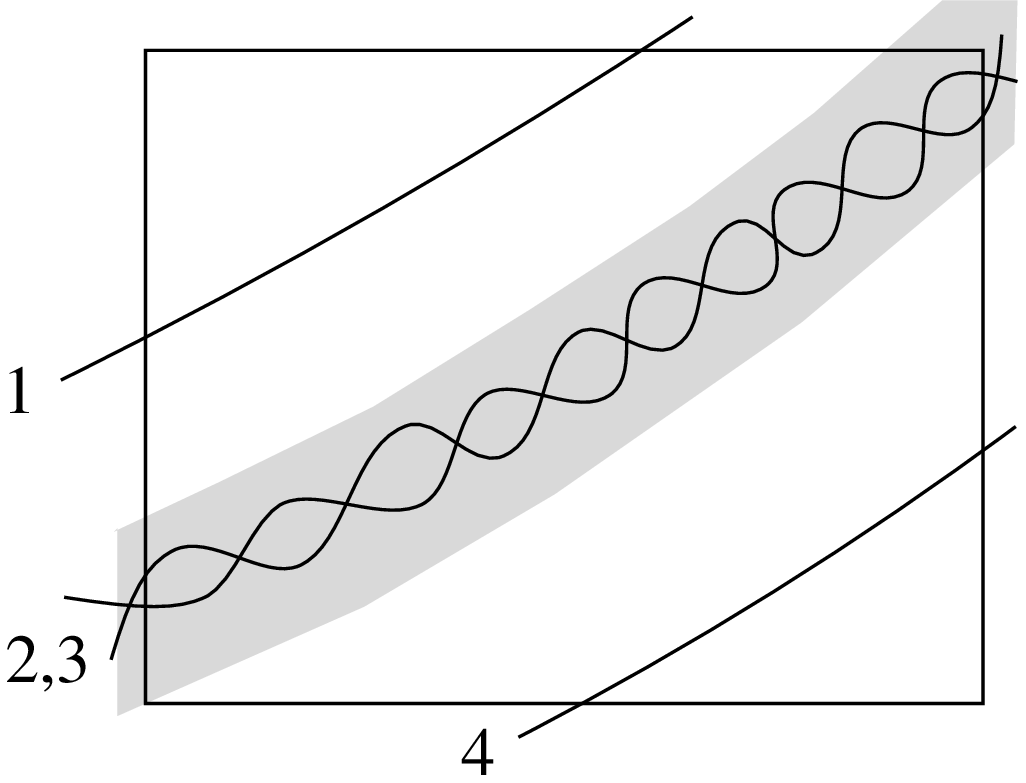}
}}\\
\sc\figurename\ 1.3&
\sc\figurename\ 1.4
\end{tabular}

The most difficult part of the proof is the estimation of $T_{jk}$
supported near the branches of $\cal Z$ having nonzero asymptotic
expansion. It is easy to see by almost orthogonality that each of
these $T_{jk}$ can be treated separately. 

A typical geometric situation is shown in Fig.~1.4:
we will have several branches of $\cal Z$ intersecting the support of
$T_{jk}$.
Actually, we show that branches which are simple, i.e.~have the 
asymptotic expansion different from asymptotic expansions of all the other branches, 
are given by functions $y_i(x)$ which are $\cin$ away from the origin 
(branches 1, 4 in Fig.~1.4). On the other hand,
multiple branches (2, 3 in Fig.~1.4) may not even be
differentiable.

In Section 5, we isolate the multiple branches by a narrow cut-off
(shaded area in Fig. ~1.4). In the rest of $\supp T_{jk}$,  we take a Whitney-type
decomposition into dyadic rectangles of the size $\gd\times L\gd$,
where $\gd$ is comparable to the distance from the rectangle to $\cal
Z$ in the anisotropic metric $|x|+L^{-1}|y|$. Here $L$ is determined
by the condition $y_i'(x)\approx L$ for the $\cin$ branches on $\supp
T_{jk}$. 

On each Whitney rectangle we have good control over $\sxy$,
which leads to an estimate for the corresponding part of $T_{jk}$. On
the other hand, we show that the rectangles with fixed $\gd$ form
an almost orthogonal family. This fact is used to obtain the optimal
norm estimate for the part of $T_{jk}$ supported away from the
multiple branches. We believe
that this argument is simpler than the inductive procedure of separating the
branches 
applied in a similar situation by Phong and Stein.

In Section 6, we deal with the part of $T_{jk}$ supported near the multiple
branches of $\cal Z$, i.e.~in the shaded region in Fig.~1.4. Here we
apply Seeger's method \cite{See1} with certain improvements possible
in our case, and obtain the norm estimate
$\gl^{-1/(N+2)}(\log\gl)^{2N/(N+2)}$, where $N$ is the multiplicity
of the branch ($N=2$ in Fig.~1.4). 
This is exactly what is claimed in \rf{main1} if $\sxy$ is completely
degenerate. 
If it is not, we show that this estimate is even better than
$\gl^{-\gd/2}$. This finishes the proof of the theorem.

\emph{Remark.} By a completely different elementary proof based on a stopping-time
argument, we can prove that for $N=2$ the estimate \rf{main1} can be
improved to the optimal $\nt\<C\gl^{-1/4}$. However, we do not know
whether a similar improvement is possible for $N\>3$.

\emph{Acknowledgements.} I wish to thank my thesis advisor E.M.~Stein
for bringing this problem to my attention, and for many stimulating discussions.

\section{Factorization of $\cin$ functions}

Recall that if $R$ is a ring, $R[t]$ and $R[[t]]$ are the rings of 
polynomials and, respectively, formal power series in $t$ with
coefficients from $R$. 
This notation can be iterated, e.g. $R[[x]][y]$ is the ring of
polynomials
in $y$ with coefficients which are elements of $R[[x]]$, and
$R[[x,y]]$ is the ring of double formal power series. 

Our aim in this section will be to prove certain factorization
formulas for $\cin$ function, which will be valid in a small
neighborhood of the origin. Since we do not care how small this
neighborhood is, it will be convenient to formulate our results for
\emph{function-germs} rather than functions. An identity involving
several function-germs is defined to be true if there exist
functions from the equivalence classes of these germs such that in the
intersection of their domains of definition the identity is true in
the usual sense.

We will make use of the following rings of germs of $\cc$-valued functions:\\  
$\bullet$ $C\(x\)$ --- continuous functions at the origin of $\rr$;\\ 
$\bullet$ $\cin\(x\)$ and $\cinxy$ --- $\cin$ functions at the origin of $\rr$
and $\rr^2$, respectively;\\
$\bullet$ $C_+\(x\)$ and $\cin_+\(x\)$ --- rings of one-sided germs; consist of
(the equivalence classes of) functions $f(x)$ 
defined in a left half-neighborhood of zero of
the form $[0,\eps)$, where $\eps>0$ can depend on $f(x)$, which are
continuous, respectively $\cin$, up to zero;

For the elements of $\cinxy$, $\cin\(x\)$, and
$\cin_+\(x\)$, we can talk about their Taylor series at the origin.
A germ whose Taylor series is zero is called \emph{flat}.
\\
$\bullet$ $A_+\(x^\gg\)$, $\gg>0$, --- the subring of 
$C_+\(x\)$ consisting of germss $f(x)$, for which there exists a
series
$\bar f (x)\in\cc[[x^\gg]]$, $\bar f(x)=\sumi_{n=0}c_n x^{n\gg}$,
such that $f(x)\sim\bar f(x)$ in the sense that for any $N$
\[
f(x)-\sum_{n=0}^N c_n x^{n\gg}=O(x^{(N+1)\gg}),\qquad x\to0.
\]
Such an $\bar f(x)$ is uniquely determined and is called the
\emph{asymptotic expansion} of $f(x)$.

The rings of germs of $\rr$-valued functions will be denoted by adding
an $\rr$ to the above notation, e.g. $\rr\cinxy$.

The main result of this section is 
the following

\begin{Prop}\label{factor}
Let $F(x,y)\in\rr\cinxy$, $\GG=\GG(F)$ be the Newton polygon of $F(x,y)$,
and assume that $\GG\ne\emptyset$. Let $\ga$ run through all 
compact edges of the boundary of $\GG$. For each edge $\ga$ joining
integer points $(A_\ga,B_\ga)$ and $(A_\ga',B_\ga')$, where $B_\ga'>B_\ga$,
put $n_\ga=B_\ga'-B_\ga$ and $\gg_\ga=(A_\ga-A_\ga')/(B_\ga'-B_\ga)$.
Let also $A$ be the $x$-coordinate of the vertical infinite edge, 
and $B$ be the $y$-coordinate of the horizontal infinite edge
of $\GG$. Then the germ $F(x,y)$ admits in the region $x,y>0$ a
factorization of the form
\begin{equation}
\label{factor1}
F(x,y)=U(x,y)\prod_{i=1}^{A}(x-X_i(y)) 
\prod_{i=1}^{B}(y-Y_i(x)) 
\prod_\ga\prod_{i=1}^{n_\ga}(y-Y_{\ga i}(x)),
\end{equation}
where\\
\textup{(1)} $U(x,y)\in\rr\cinxy$, $U(0,0)\ne 0$, \\  
\textup{(2)} all $X_i(x), Y_i(x)\in C_+\(x\)$, and
$X_i(x), Y_i(x)=O(x^N)$ as $x\to 0$ for any $N>0$,\\ 
\textup{(3)} all $Y_{\ga i}(x)\in A_+\(x^{1/n!}\)$ for $n=
B+\sum_\ga n_\ga$ with asymptotic expansions of the form $Y_{\ga i}(x)=
c_{\ga i} x^{\gg_\ga}+\ldots$ as $x\to 0$, where $c_{\ga i}\ne 0$,\\
\textup{(4)} if $Y(x)$ is any of the functions $Y_{\ga i}(x)$, and if 
$f(x,y)=\prod (y-Y_{\ga i}(x))$ is the product over all $i$ such that 
$Y_{\ga i}(x)$ has exactly the same asymptotic expansion as $Y(x)$,
then $f(x^{n!},y)\in \cin_+\(x\)[y]$,\\
\textup{(5)} if in \textup{(4)} we additionally assume that the
asymptotic expansion of $Y(x)$ is real, then $f(x,y)$ is also real.
\end{Prop}

The proof relies on the following result, which is well known in
the theory of plane algebraic curves as the \emph{Puiseux theorem} (see
\cite{Wal}, p.~98ff, or \cite{Bour}, A.V.150).

\begin{Lem}
\label{puiseux} Let $\bar F(x,y)\in\cc[[x]][y]$ be of the form
\[
\bar F(x,y)=y^n+
\bar c_{n-1}(x) y^{n-1}+\ldots+\bar c_0(x),\qquad \bar c_i(x)\in\cc[[x]],
\]
where the zeroth order terms of all $c_i(x)$ vanish. Let $\ga$, $n_\ga$, $\gg_\ga$, $B$
be
defined via the Newton polygon $\GG=\GG(\bar F)$ in the same way as in
Proposition {\rm\ref{factor}}. Then there exists a factorization
\begin{equation}
\label{puiseux1}
\bar F(x,y)=y^{B}\prod_\ga\prod_{i=1}^{n_\ga}(y-\bar Y_{\ga i}(x)),
\end{equation}
where the series $\bar Y_{\ga i}(x)\in\cc[[x^{1/n!}]]$ are of the
form $\bar Y_{\ga i}(x)=c_{\ga i} x^{\gg_\ga}+\ldots$ with $c_{\ga
i}\ne 0$.
\end{Lem}

The following lemma will be used to get factorizations of
function-germs from
factorizations of their formal Taylor series.

\begin{Lem} 
\label{poly}
Let $P(x,y)\in\cin\(x\)[y]$ be of the form
\[
P(x,y)=y^n+
 c_{n-1}(x) y^{n-1}+\ldots+c_0(x),\qquad c_i(x)\in\cin\(x\),
\]
where $c_i(0)=0$ for all $i$. Let $\bar P(x,y)\in\cc[[x]][y]$ 
be the formal Taylor series of $P(x,y)$ at the origin. 
Let $\bar Y(x)\in\cc[[x]]$ be a root of multiplicity $m$, $1\<m\<n$,
of $\bar P(x,y)$ considered as a polynomial in $y$, i.e.
\[
\bar P(x,\bar Y(x))=\ldots=\bar P_y^{(m-1)}(x, \bar Y(x))=0,
\qquad \bar P_y^{(m)}(x, \bar Y(x))\ne 0
\]
as elements of $\cc[[x]]$. Then there exist $m$ function-germs 
$Y_1(x),\ldots,Y_m(x)\in C\(x\)$ such that\\
\textup{(1)} all $Y_i(x)\sim \bar Y(x)$ as $x\to 0$,\\
\textup{(2)} all $P(x,Y_i(x))=0$ for $x>0$,\\
\textup{(3)} $\prod_{i=1}^m(y-Y_i(x))\in\cin\(x\)[y]$,\\
\textup{(4)} if we additionally assume that $P(x,y)$ and $\bar Y(x)$
are real, then $\prod_{i=1}^m(y-Y_i(x))$ is also real.
\end{Lem}
\begin{proof} Let $\tilde Y(x)$ be a $\cin$ function with the formal
Taylor series $\bar Y(x)$, supplied by E.~Borel's theorem. Denote 
\[
\delta_i(x)=\frac 1{i!}P_y^{(i)}(x,\tilde Y(x)),\qquad i=0,\ldots,n.
\]
Let the (nonzero by assumption) series $\bar P^{(m)}_y(x,\bar Y(x))$
starts with a term $cx^s$, $c\ne 0$, $s\in \zp$. 
Then we have $\gd_m(x)=cx^s+o(x^s)$ as $x\to 0$. On the other
hand, it is clear that $\gd_0(x),\ldots,\gd_{m-1}(x)$ are flat.

We will be looking for $Y_i(x)$ of the form 
\[
Y(x)=\tilde Y(x)+\delta_m(x)\ga(x),
\]
where $\ga(x)$ is an unknown continuous $\cc$-valued function-germ such that
$\ga(x)=O(x^N)$ as $x\to0$ for any $N>0$.

By Taylor's formula, the equation $P(x,Y(x))=0$ can be written as
\begin{equation}
\sum_{i=0}^n \gd_i(x)[\gd_m(x)\ga(x)]^i=0.
\end{equation}
For small $x$, this is equivalent to the equation
$w(x,\ga(x))=0$ for the function $w(x,z)$ given by
\[
w(x,z)=\sum_{i=0}^n \gd_i(x)[\gd_m(x)]^{i-m-1} z^i.
\]
Note that if $f,g\in\cin$, and $f$ is flat at the origin, while $g$ is
not flat, then $f/g$ is $\cin$ near the origin and is flat. It follows that
$w(x,z)\in\cin\(x\)[z]$. 

 On the complex circle $|z|=x^N$,
$N>0$, the term $z^m$ will dominate the other terms in $w(x,z)$ if $x$ is
sufficiently small. By Rouche's theorem it follows that the equation
$w(x,z)=0$ has for small fixed $x$ exactly $m$ roots in the disc
$|z|<x^N$, which we denote
$\ga_i(x)$, $i=1,\ldots,m$. We can arrange so that $\ga_i(x)$ are
continuous in $x$, and the previous argument shows that
$\ga_i(x)=O(x^N)$ for any $N>0$.

We now prove (3). Since the functions $\ga_i(x)$ 
enter the product
in (3) in a symmetric way, it is sufficient to prove that the
elementary symmetric polynomials $s_0,\ldots, s_m$ in $\ga_i(x)$ are in
$\cin\(x\)$. By the Newton relations (see
\cite{Bour}, A.IV.70), it is sufficient to prove the same for the
functions
\[
p_k(x)=\sum_{i=1}^m [\ga_i(x)]^k,\qquad k=1,\ldots,m.
\]
However, by Cauchy's formula we have that for small $x$
\[
p_k(x)=\frac1{2\pi i}\oint_{|z|=\eps} \frac {z^k w'_z(x,z)}{w(x,z)}\,dz,
\]
from where it is clear that $p_k(x)\in\cin\(x\)$.

To prove (4), we notice that under the additional assumption made we
can take $\tilde Y(x)$ to be real. Then $w(x,z)\in\rr\cin\(x\)[z]$,
and therefore non-real roots $\ga_i(x)$ will appear in conjugate
pairs. Then all $p_k(x)$ will be real, which implies (4). 
\end{proof}

Now we can prove 

\begin{Lem} 
\label{pfactor}
Proposition {\rm\ref{factor}} is true if
$F(x,y)\in\rr\cin\(x\)[y]$.
\end{Lem} 
\begin{proof}
By Lemma \ref{puiseux}, the Taylor series $\bar F(x,y)\in\rr[[x]][y]$ 
of $F(x,y)$ has a factorization \rf{puiseux1}. Consider the function 
$P(x,y)=F(x^{n!},y)$. Its Taylor series 
has the form $\bar P(x,y)=\bar F(x^{n!},y)$, and so factorizes as
\[
 \bar P(x,y)=y^{B}\prod_\ga\prod_{i=1}^{n_\ga}(y-\bar Y_{\ga i}(x^{n!})).
\]
Let $\bar Y(x)$ be one of the series $\bar Y_{\ga i}(x^{n!})\in\cc[[x]]$,
and assume that among all the $\bar Y_{\ga i}(x^{n!})$ there are exactly $m$
series coinciding with $\bar Y(x)$. Then $y=\bar Y(x)$ is a root of
multiplicity $m$ of the polynomial $\bar P(x,y)\in\rr[[x]][y]$, and by
Lemma \ref{poly} we conclude that there exist $m$ functions $Y_i(x)\in
C\(x\)$, $i=1,\ldots,m$, such that (1)--(3) from the formulation of
the lemma are true. 

In view of (3), we can divide $P(x,y)$ by
$\prod(y-Y_i(x))$, and the result is again a polynomial $\tilde
P(x,y)$ from $\cin\(x\)[y]$. The Taylor polynomial of $\tilde
P(x,y)$ will be $\bar P(x,y)$ divided by $(y-\bar Y(x))^m$. Now we can
apply Lemma \ref{poly} to $\tilde P(x,y)$ choosing a different $\bar
Y(x)$ etc. 

By repeating this operation several times, we get a complete factorization of $P(x,y)$. The
required factorization of $F(x,y)$ is then obtained by the inverse
substitution $x\mapsto x^{1/n!}$. The property (5) is ensured by 
splitting off all real series $\bar Y(x)$ before non-real ones in the
above argument.  
\end{proof}

Proposition \ref{factor} will be reduced to Lemma \ref{pfactor}
by means of the following \emph{Malgrange preparation theorem} (see
\cite{GG}, p.~95).

\begin{Lem}
\label{malgrange}
Let $F(x,y)\in\rr\cinxy$, and assume that $F(0,y)$ is not flat, so that
$F(0,y)=cy^n+o(y^n)$, $y\to0$, for some $n\in\zp$ and $c\ne0$. Then
there is a factorization 
\[
F(x,y)=U(x,y)P(x,y),
\]
where \\
\textup{(1)} $U(x,y)\in\rr\cinxy$, $U(0,0)\ne 0$,\\
\textup{(2)} $P(x,y)\in\rr\cin\(x\)[y]$ is of the form
\[
P(x,y)=y^n+c_{n-1}(x)y^{n-1}+\ldots+c_0(x),
\]
where all $c_i(x)\in\rr\cin\(x\)$, $c_i(0)=0$.
\end{Lem}

\begin{proof}[Proof of Proposition {\rm\ref{factor}}.]
Notice that the Newton polygon is invariant with \mbox{respect} to
multiplication by a nonzero $\cin$ function (see Phong and Stein
\cite{PS}, p.~112). Therefore, for the functions $F(x,y)$ such that $F(0,y)$
is not flat (which is equivalent to having $A=0$) the proposition
follows immediately from Lemmas \ref{malgrange} and \ref{pfactor}.

Assume now that $A>0$. In this case we must somehow separate the roots
infinitely tangent to the $y$-axis. This can be done as follows. Since
$F(x,y)$ is not flat at the origin, there exists a rotated orthogonal
system of coordinates $(x',y')$ such that the restriction of $F$ to
the $y'$-axis is not flat. So we can apply Lemma \ref{malgrange} to
$F$ written in coordinates $(x',y')$. Let $P(x',y')$ be the arizing
polynomial.    

If $y'=ax'$ is the equation of the old $y$-axis in the new
coordinates, then $y'=ax'$ will be a root of multiplicity $A$ of $\bar
P(x',y')\in\rr[[x']][y']$. So we can apply Lemma \ref{poly} and obtain
$A$ roots $y'=Y_i(x')$, $i=1,\ldots,A,$ of $P(x',y')=0$, such that
$Y_i(x')\sim ax'$.

Moreover, by Lemma \ref{poly} (3),(4) we will have that
$Q(x',y')=\prod(y'-Y_i(x'))$ is in $\rr\cin\(x'\)[y']$. So we can
divide $P(x',y')$ by $Q(x',y')$, and the quotient will be a $\cin$
function, which is no longer flat on the old $y$-axis.

Let $\tilde F(x,y)$ be this last quotient written in the old system of
coordinates. Then $\GG(\tilde F)$ is just $\GG(F)$ shifted $A$
units to the left. So we can factorize $\tilde F(x,y)$ as in the case
$A=0$ described above. 

It remains to get a factorization of $Q(x',y')$ in the old
coordinates. It is clear that the Taylor series of $Q$ written in
the coordinates $(x,y)$ consists of one term $cx^A$. Interchanging
the roles of $x$ and $y$ brings us back to the case $A=0$, and the
required factorization of the form $\prod(x-X_i(y))$ can be obtained
as described above.
\end{proof}


\section{Dyadic decomposition of $T$. Estimates away from $\cal Z$.}

We decompose the operator $T$ as
\[
T=\sum_{\pm}\sum_{j,k} T_{jk}^{\pm\pm},
\] 
where $T_{jk}^{++}$ is defined as 
\[
T_{jk}^{++}f(x)
=\int_{-\infty}^{\infty}e^{i\gl S(x,y)}\chi_j(x)\chi_k(y)\chi(x,y)f(y)\,dy. 
\] 
Here $\sum_j\chi_j(x)=1$ is a smooth dyadic partition of unity on
$\rr^+$,
so that the kernel of $T_{jk}^{++}$ is supported on the rectangle
$R_{jk}=[2^{-j-1},2^{-j+1}]\times[2^{-k-1},2^{-k+1}]$. Three other
$\pm$
combinations refer to the quadrants defined by specific signs of $x$
and $y$. We restrict ourselves with the positive quadrant, the other
ones
being exactly similar, and denote $T_{jk}^{++}$ by simply $T_{jk}$.

Denote $F(x,y)=\sxy(x,y)$. By Lemma \ref{factor}, there is 
a neighborhood of the origin $V$ such that in $V\cap \rrp^2$ 
there exists a factorization of the form \rf{factor1}. 
We assume that $\supp \chi\subset V$. The singular variety 
$\cal Z=\{(x,y)\in V:F(x,y)=0\}$ now splits into branches corresponding
to the factors in the RHS of \rf{factor1} (see Fig.~1.1). Note, however, that some
of these branches may contain an imaginary component.

Let $R_{jk}^*$ denote the double of $R_{jk}$.
We fix a large constant $D$ such that if the pair $(j,k)$
satisfies the condition $\min_\ga|k-j\gg_\ga|\>D$, then $y-c_{\ga
i}x^{\gg_\ga}\ne 0$ on $R_{jk}^*$ for all $c_{\ga i}x^{\gg_\ga}$
occurring
as the lowest order terms of the asymptotic expansions of $Y_{\ga i}(x)$
in Proposition \ref{factor}(3). 

Let us number the compact edges $\ga$ of the boundary of the Newton polygon $\GG(F)$ 
from left to right, so that $\gg_1<\gg_2<\ldots<\gg_{\ga_0}$, where
$\ga_0$ 
is the total number of compact edges. Also put $\gg'=\gg_1/2$ if
$A>0$, $\gg'=0$ otherwise; $\gg''=2\gg_{\ga_0}$ if
$B>0$, $\gg''=\infty$ otherwise. 

Consider the operators
\begin{equation}
\label{talpha}
T_\nu=\sum_{j\gg_\nu\ll k\ll
j\gg_{\nu+1}} T_{jk},\qquad  1\<\nu\<\ga_0-1,
\end{equation} 
\[
T'=\sum_{j\gg'\ll k\ll
j\gg_{1}} T_{jk},\qquad
T''=\sum_{j\gg_{\ga_0}\ll k\ll
j\gg''} T_{jk},
\]
where $a\ll b$ stands for $a\<b-D$. 
These operators constitute the part of $T$ supported relatively
far away from $\cal Z$ (see Fig.~1.2). 
In this section, we prove that $\|T_\nu\|\<C\gl^{-\gd/2}$ for each $\nu$. 
The reader will believe us that with minor modifications the
argument given below will also produce the same estimate for $T'$, $T''$.
\begin{Lem} \label{operatorcorput} 
Let $T$ be an oscillatory integral operator of the form
\rf{operator}. Assume that\\
\tup{1} $\chi(x,y)$ is supported in a
rectangle $R$ of size $\gd_x\times\gd_y$,\\
\tup{2} $|\D_y^n\chi|\<C/\gd_y^{n}$  in  $R$ for $n=0,1,2$,\\
\tup{3} $|\sxy|\>\mu>0$  in  $R$,\\
\tup{4} $|\D_y^n\sxy|\<C\mu/\gd_y^n$  in  $R$ for $n=1,2$.\\
Then $\nt\<const(\gl\mu)^{-1/2}$ with $const$ depending only on $C$.
\end{Lem}

This is a variant of the \emph{Operator van der Corput lemma} of Phong
and
Stein \cite{PS}. The lemma is proved by a standard $TT^*$
argument. The assumptions made are enough to show, integrating
by parts twice, that the kernel of $TT^*$ has the bound
\[
K(x_1,x_2)\<C_A\frac{\gd_y}{1+\gl^2\mu^2\gd_y^2|x_1-x_2|^2},
\]
which implies the necessary norm estimate. We omit the details.

We will need the following partial case of the Schur test (see
\cite{Halmos}, Theorem 5.2).

\begin{Lem}
\label{schur} 
Let $T$ be an integral operator on $L^2(\rr)$ with kernel $K(x,y)$,
\[
Tf(x)=\int_{-\infty}^\infty K(x,y)f(y)\,dy.
\]
Assume that the quantities
\begin{align*}
M_1=\sup_y\int|K(x,y)|\,dx,\qquad
M_2=\sup_x\int|K(x,y)|\,dy
\end{align*}
are finite. Then $T$ is bounded with $\nt\<(M_1 M_2)^{1/2}$.
\end{Lem}
    

\begin{Cor}
\label{size}
 \tup{Phong and Stein \cite{PSmodels}, Lemma 1.6}
Let $T$ be an integral operator with kernel $K(x,y)$, and assume that\\
\tup{1} $|K(x,y)|\<1$,\\
\tup{2} for each $y$, $K(x,y)$ is supported in an $x$-set of measure
$\<\gd_x$,\\
\tup{3} for each $x$, $K(x,y)$ is supported in a $y$-set of measure
$\<\gd_y$.\\
Then $\nt\<(\gd_x\gd_y)^{1/2}$.
\end{Cor}  

It is natural to call the latter bound on $\|T\|$ the \emph{size
estimate}, as opposed to the \emph{oscillatory estimate}
supplied by Lemma \ref{operatorcorput}.

In the rest of the paper we will use the notation $a\lesssim b$
to
mean $a\<Cb$, and
$a\approx b$ to mean $C^{-1}b\<a\<Cb$, where $C>0$ is an unimportant
constant, which is supposed to be independent of $j,k,\gl$.

Take an operator $T_{jk}$ entering
the RHS of \rf{talpha}. We may reduce $V$ if necessary so that on the 
part of $\cal Z$ inside $V$ the functions $Y_{\ga i}$ do not differ much
from the first terms of their asymptotic expansions. Assume that 
$T_{jk}$ is nonzero, which means that $R_{jk}\cap V\ne\emptyset$.
Then it is clear from the definition of $D$ that the
factors in the RHS of \rf{factor1} can estimated as follows for 
$(x,y)\in R_{jk}$:
\begin{equation}
\label{boundsaway}
\begin{array}{l}
|x-X_i(y)| \approx 2^{-j},\\
|y-Y_i(x)| \approx 2^{-k},\\
|y-Y_{\ga i}(x)| \approx 
\begin{cases}
2^{-k}, &\ga>\nu,\\
2^{-j\gg_\ga},&\ga\<\nu.
\end{cases}
\end{array}
\end{equation}
Therefore it follows from \rf{factor1} that on $R_{jk}$
\begin{equation}
\label{hessian}
|F|\approx 2^{-jA}2^{-kB}\prod_{\ga>\nu} 2^{-kn_\ga}
\prod_{\ga\<\nu}
2^{-j\gg_\ga n_\ga}=:\mu.
\end{equation}
The numbers $\gg_\ga$, $n_\ga$ can be found from the Newton polygon 
$\GG(F)$ as described in Proposition \ref{factor}. By using this 
information, we find that
\[
\mu=2^{-jA_\nu-kB_\nu}.
\]
We claim that on $R_{jk}$
\begin{equation}
\label{derivatives}
|\D_y^n F|\lesssim \mu2^{kn},\qquad n=1,2.
\end{equation}
Indeed, when differentiating the RHS of \rf{factor1} in $y$, the derivative
can fall one either $U(x,y)$, or $\prod(x-X_i(y))$, 
or one of the 
remaining terms. In the first case, we simply get a bounded factor. In the
second
case, we get an even better factor of $O(2^{-kN})$ for any $N>0$, since 
the product in question is a $\cin$ function whose Taylor series
at the origin is $x^{A}$. Finally, in the third case we get a
factor of the form $(y-Y_i(x))^{-1}$ or $(y-Y_{\ga i}(x))^{-1}$, which
is $O(2^k)$ in view of \rf{boundsaway}. This argument works equally
well
for the second derivative, giving \rf{derivatives}.

The rectangle $R_{jk}$ is of size $\gd_x\times\gd_y$ with $\gd_x\approx
2^{-j}$,
$\gd_y\approx2^{-k}$. So the conditions of Lemma \ref{operatorcorput}
are satisfied, and we obtain the oscillatory estimate 
\begin{equation}
\label{oscillatory}
\|T_{jk}\|\lesssim \gl^{-1/2}2^{(jA_\nu+kB_\nu)/2}.
\end{equation}
On the other hand, the size estimate following from Corollary \ref{size} is
\begin{equation}
\label{sizeest}
\|T_{jk}\|\lesssim 2^{(j+k)/2}.
\end{equation}
The required bound for $T_\nu$ can now be derived from the last two
estimates by summation, taking into account possible orthogonality
relations between different $T_{jk}$. We use the summation procedure
described in detail in Phong and Stein \cite{PS}, pp.~120--122. Let us
briefly recall what that procedure looks like. 
 
We have three different cases: $A_\nu>B_\nu$, $A_\nu<B_\nu$, and
$A_\nu=B_\nu$.
If $A_\nu>B_\nu$, we put $k=[\gg_\nu j]+r$, $r\>-D$, 
substitute this into \rf{oscillatory} and \rf{sizeest}, and take 
the geometric mean of the two estimates killing the
$j$-factor.
The result is (\cite{PS}, eq.~(4.35))
\[
\|T_{jk}\|\lesssim\gl^{-\gd_\nu/2}2^{-(1-\gd_\nu-B_\nu\gd_\nu)r/2},
\]
where
\begin{equation}
\label{deltanu}
\gd_\nu=\frac{1+\gg_\nu}{1+A_\nu+(1+B_\nu)\gg_\nu}.
\end{equation}
For a fixed $r$, the same estimate is true by almost orthogonality 
for the sum of $T_{jk}$ over $(j,k)$ satisfying $k=[\gg_\nu j]+r$,
since operators $T_{jk}$ and $T_{j'k'}$ in such a sum have
disjoint $x$- and $y$-supports for $|j-j'|$ larger than a fixed
constant:
\begin{equation}
\label{fixedr}
\biggl\|\sum _{k=[\gg_\nu j]+r}T_{jk}\biggr\|\lesssim\gl^{-\gd_\nu/2}2^{-(1-\gd_\nu-B_\nu\gd_\nu)r/2},
\end{equation}
In the case under consideration
\[
1-\gd_\nu-B_\nu\gd_\nu=\frac{A_\nu-B_\nu}{1+A_\nu+(1+B_\nu)\gg_\nu}>0,
\]
so we can sum \rf{fixedr} in $r$, and obtain 
\begin{equation*}
\|T_\nu\|\<\sumi_{r=-D}\biggl\|\sum _{k=[\gg_\nu
j]+r}T_{jk}\biggr\|\lesssim
\gl^{-\gd_\nu/2}.
\end{equation*}
Let $t_\nu$ be the parameter of intersection of the line $n_1=n_2=t$
with the line containing edge $\nu$ of $\GG(F)$. Then a simple
calculations shows that $1/\gd_\nu=1+t_\nu$. Therefore it is clear
that $\gd_\nu\>\gd$, and the previous estimate implies 
$\|T_\nu\|\lesssim\gl^{-\gd/2}$, as desired.

Further, the case $A_\nu<B_\nu$ is analogous to the case
$A_\nu>B_\nu$ and in fact reduces to it by interchanging the roles
of $j$ and $k$.

Finally, in the case $A_\nu=B_\nu$ we first sum the estimates 
\rf{oscillatory} and \rf{sizeest} along the diagonals $j+k=i$. The 
corresponding sums are again almost orthogonal, and we get
\[
\biggl\|\sum _{j+k=i}T_{jk}\biggr\|\lesssim \min(2^{-i/2},\gl^{-1/2}2^{iA_\nu/2}),
\]
whence
\[
\|T_\nu\|\<\sumi_{i=0}\biggl\|\sum _{j+k=i}T_{jk}\biggr\|\lesssim
\gl^{-1/(2+2A_\nu)}.
\]
Since in the considered case $1/(1+A_\nu)=\gd_\nu$, we again recover
the required estimate.

The treatment of the operators $T_\nu$ is now complete.

\section{Estimates near the coordinate axes.}

In this section, we deal with the operators
\[
T_y=\sum_{k\ll\gg'j} T_{jk},\qquad T_x=\sum_{\gg''j\ll k} T_{jk}.
\]
These two constitute the part of $T$ supported near the branches of
$\cal Z$
which are infinitely
tangent to the coordinate axes (see Fig.~1.3). We will prove the estimate
$\|T_x\|\lesssim
\gl^{-\gd/2}$. The same estimate will be true for $T_y$, since 
taking the adjoint of $T$ brings $T_y$ to the form of $T_x$. 
We may of course assume $B\>1$, since otherwise $\gg''=\infty$ and $T_x=0$.

We represent $T_x$ as
\[
T_x=\sum_j T_j,\qquad T_j=\sum_{\gg''j\ll k} T_{jk},
\]
and claim that\\
(1) $\|T_j\|\lesssim\gl^{-\gd/2},$\\
(2) $\|T_j^*T_{j'}\|=0$ for $|j-j'|\>2$,\\
(3) $\|T_j T_{j'}^*\|\lesssim\gl^{-\gd/2}2^{-\eps|j-j'|}$ for some
$\eps>0$.\\
If we prove all these, the estimate $\|T_x\|\lesssim
\gl^{-\gd/2}$ will follow from the Cotlar--Stein lemma.

We have 
\[
T_{j}f(x)
=\int_{-\infty}^{\infty}e^{i\gl
S(x,y)}\chi_j(x)\tilde\chi_j(y)\chi(x,y)f(y)\,dy,  
\]
where $\tilde \chi_j=\sum_{\gg''j\ll k}\chi_k$, so that the support
of $\tilde\chi$ is contained in $[0,C2^{-\gg'' j}]$. The property (2) is
obvious. Further, the operator $T_jT_{j'}^*$ has the kernel
\[
K(x_1,x_2)=\chi_j(x_1)\chi_{j'}(x_2)\int e^{i\gl[S(x_1,y)-S(x_2,y)]}
\tilde\chi_j(y)\tilde\chi_{j'}(y)\chi(x_1,y)\chi(x_2,y)\,dy.
\] 
We want to estimate this by the following variant of the standard van der
Corput lemma (see \cite{STEIN}, Corollary on p.~334).

\begin{Lem} 
\label{vandercorput}
Let $k$ be a positive integer, $\GF\in C^{k}[a,b]$, $\GP\in C^1[a,b]$, and assume that
$\GF^{(k)}\>\mu>0$ on $[a,b]$. If $k=1$, assume
additionally that $\GF'$ is monotonic on $[a,b]$. Then
\[
\biggl|\int_a^b e^{i\gl\GF(y)}\GP(y)\,dy\biggr|\lesssim
(\gl\mu)^{-1/k}\biggl(|\GP(b)|+\int_a^b |\GP'|\biggr).
\]
\end{Lem}

Assume that $j'\>j$. We apply this lemma with $[a,b]=[0,C2^{-\gg''j}]$, 
 $k=B+1\>2$,
\begin{align*}
\GF(y)&=S(x_1,y)-S(x_2,y),\\
\GP(y)&=\tilde\chi_j(y)\tilde\chi_{j'}(y)\chi(x_1,y)\chi(x_2,y).
\end{align*} 
It is clear that $|\GP(b)|+\int_a^b |\GP'|\lesssim 1$. Further (recall
that we denoted $\sxy=F$), 
\[
\GP^{(B+1)}(y)=\D_y^{B+1} S(x_1,y)-\D_y^{B+1} 
S(x_2,y)=\int_{x_2}^{x_1}\D_y^{B}
F(x,y)\,dx.
\]
Of all the terms arising when we differentiate \rf{factor1} 
$B$ times in $y$, the term in which all derivatives fall on  
$\prod(y-Y_j(x))$ will dominate on the support of $T_x$ after 
a possible reduction of $V$. It follows
that on the support of $T_x$ 
\[ 
|\D_y^{B}F|\approx x^{A+\sum_\ga n_\ga\gg_\ga}=x^{A'},
\]
where we denoted $A'=A_{\ga_0}$. Note that $(A',B)$ is the common vertex
of the horizontal infinite edge of $\GG(F)$ and its last compact edge $\ga_0$.

By the previous remarks,
\[
|\GP^{(B+1)}(y)|\gtrsim\bigl|x_1^{A'+1}-x_2^{A'+1}\bigr|.
\]
In the case $j'=j$ we have $
\bigl|x_1^{A'+1}-x_2^{A'+1}\bigr|\approx 2^{-jA'}
|x_1-x_2|$ on the support of $K(x_1,x_2)$, so Lemma \ref{vandercorput} 
gives
\[
|K(x_1,x_2)|\<2^{jA'/(
 B+1)}(\gl
|x_1-x_2|)^{-1/(
 B+1)}.
\]
So by Lemma \ref{schur},
\begin{align*}
\|T_j T_j^*\|&\lesssim 2^{jA'/(
 B+1)}\int_0^{2^{-j}}(\gl t)^{-1/(
 B+1)}\,dt\\
&\lesssim \gl^{-1/(B+1)} 2^{j(A'-B)/(B+1)}.
\end{align*}
This of course implies the estimate
\begin{equation}
\label{est1}
\|T_j\|\lesssim  \gl^{-1/(2B+2)} 2^{j(A'-B)/(2B+2)}.
\end{equation}
Another estimate is supplied by Corollary \ref{size}:
\begin{equation}
\label{est2}
\|T_j\|\lesssim 2^{-j(1+\gg'')/2}\<2^{-j(1+\gg_{\ga_0})/2}.
\end{equation}
As the reader may check, 
taking the geometrical mean of these two bounds which
kills the $j$-factor gives exactly $\|T_j\|\lesssim
\gl^{-\gd_{\ga_0}/2}$,
with $\gd_{\ga_0}$ defined as in \rf{deltanu}. This implies (1) since
all $\gd_\nu\>\gd$.

In proving (3), we may assume $j'\>j+2$. Then 
$
\bigl|x_1^{A'+1}-x_2^{A'+1}\bigr|\approx 2^{-j(A'+1)}$ 
on the support of $K(x_1,x_2)$, whence by Lemma \ref{vandercorput}
\[
|K(x_1,x_2)|\lesssim \gl^{-1/(B+1)}2^{j(A'+1)/(B+1)}=:M.
\]
The support of $K(x_1,x_2)$ is contained in the rectangle of size
$\approx 2^{-j}\times2^{-j'}$. Now Corollary \ref{size} 
gives the bound
\[
\|T_jT_{j'}^*\|\lesssim M2^{-(j+j')/2}=\gl^{-1/(B+1)}2^{j(A'-B)/(B+1)}
2^{-\Delta j/2},
\]
where we denoted $\Delta j=j'-j$. By multiplying the estimates
\rf{est2}
for $T_j$ and $T_{j'}$, we get another bound:
\[
\|T_jT_{j'}^*\|\lesssim 2^{-j(1+\gg_{\ga_0})}
2^{-\Delta j(1+\gg_{\ga_0})/2}.
\] 
These two bounds have the form of \rf{est1} and \rf{est2} squared, but
with
an additional factor exponentially decreasing in $\Delta j$. Therefore
it is clear that this time taking the geometric mean killing the
$j$-factor
will give 
\[
\|T_jT_{j'}^*\|\lesssim \gl^{-\gd_{\ga_0}}
2^{-\eps\Delta j} 
\]  
for some $\eps>0$. This implies (3) and concludes the treatment of $T_x$.

\section{Estimates near $\cal Z$.}

We still have to estimate the part of $T$ supported near the branches
of $Z$ which are not infinitely tangent to the coordinate axes. This
part
is the sum over $\nu=1,\ldots,\ga_0$ of the operators
\begin{equation}
\label{tnu}
T^\nu=\sum_{\gg_\nu j-D<k<\gg_\nu j+D} T_{jk}.
\end{equation}
Notice that the sum in \rf{tnu} is almost orthogonal, since the $x$- and
$y$-supports of $T_{jk}$ and $T_{j'k'}$ are disjoint for $|j-j'|$
larger than a fixed constant. Therefore it suffices to prove the
estimate $\|T_{jk}\|\lesssim\gl^{-\gd/2}$ for each $T_{jk}$ from 
the RHS of \rf{tnu}.

Fix such a $T_{jk}$. Analogously to \rf{hessian}, on $R_{jk}$
\begin{align}
\label{fonrjk}
|F|&\approx 2^{-jA}2^{-j\gg_\nu B}\prod_{\ga>\nu} 2^{-j\gg_\nu n_\ga}
\prod_{\ga<\nu}
2^{-j\gg_\ga n_\ga}\prod_{i=1}^{n_\nu}|y-Y_{\nu i}(x)| \\
&=2^{-j(\gg_\nu B_\nu+A_\nu-\gg_\nu n_\nu)}\prod_{i=1}^{n_\nu}|y-Y_{\nu
i}(x)| \notag\\
\notag&=2^{-j(\gg_\nu B_\nu+A_\nu-\gg_\nu n_\nu')}\prod_{i=1}^{n_\nu'}|y-Y_{\nu
i}(x)|,
\end{align}
where we ordered $Y_{\nu i}$ so that for $n_\nu'<i\le n_\nu$ we have 
$\Re c_{\nu i}=0$ in $Y_{\nu i}=c_{\nu i}x^{\gg_\nu}+\ldots$.

Let us quickly dispose of the case $n_\nu'=0$, in which we can apply
Lemma \ref{operatorcorput} (the condition (4) is easily checked) and
Corollary \ref{size} to get the oscillatory and size estimates
\begin{align*}
&\|T_{jk}\|\lesssim\gl^{-1/2}2^{j(\gg_\nu B_\nu+A_\nu)/2},\\
&\|T_{jk}\|\lesssim 2^{-j(1+\gg_\nu)/2}.
\end{align*}
Now by taking the geometric mean killing the $j$-factor, we obtain the 
required estimate $\|T_{jk}\|\lesssim\gl^{-\gd_\nu/2}\<\gl^{-\gd/2}$.  

Now assume that $n_\nu'>0$. Denote $r_i(x)=\Re Y_{\nu i}(x)$, and let
$\bar r _i(x)\in\rr[[x^{1/n!}]]$ be the asymptotic expansion of
$r_i(x)$ at zero. By E.~Borel's theorem, we can find real functions
$f_i(x)$ such that $f_i(x^{n!})\in\cin$ and $f_i(x)\sim \bar r_i(x)$
as $x\to 0$. Moreover, there is one case when we may and will take
simply $f_i(x)=Y_{\nu i}(x)$. Namely, by Proposition
\ref{factor} (4),(5)  
this is possible if the series 
$\bar Y_{\nu i}(x)$ is real and different from any other $\bar
Y_{\nu i'}(x)$.

Let $W$ be the union of the graphs of $f_i(x)$ inside $R_{jk}$:
\[
W=\bigcup_{i=1}^{n_\nu'}\bigl\{(x,y)\in R_{jk}\bigl|
y=f_i(x)\bigr\}.
\]
It is not difficult to see that on $R_{jk}$
\[
f_i'(x)\approx x^{\gg_{\nu}-1}\approx2^{-j(\gg_\nu-1)}=:L.
\]
This suggests to consider a Whitney-type decomposition of
$R_{jk}\backslash W$ away from $W$ into rectangles of the size $2^{-m}\times
L2^{-m}$. The easiest way to do this is to dilate the set
$R_{jk}\backslash W$ along the $y$-axis $L^{-1}$ times, take the
standard
Whitney decomposition into the dyadic squares away from (the dilation
of) $W$, and
contract
everything to the original scale. As a result, we get a covering
\[
R_{jk}\backslash W\subset \bigcup R_l,\qquad R_l\cap R_{jk}\ne\emptyset,
\]
where $R_l$ are rectangles of the size $2^{-m_l}\times
L2^{-m_l}$, $m_l\in\zp$, such that the distance from $R_l$ to $W$ in
the
anisotropic norm $|x|+L^{-1}|y|$ is of the order $2^{-m_l}$.

We claim that the rectangles $R_l$ of \emph{fixed} size form an
almost orthogonal family, i.e.~that for each $R_l$ the number of
rectangles $R_{l'}$ with $m_{l'}=m_l$ such that either the $x$- or the
$y$-projections of $R_l$ and $R_{l'}$ intersect is bounded by a fixed
constant independent of $l$. 

\begin{center}
\scalebox{0.5}{
\rotatebox{-90}{
\includegraphics{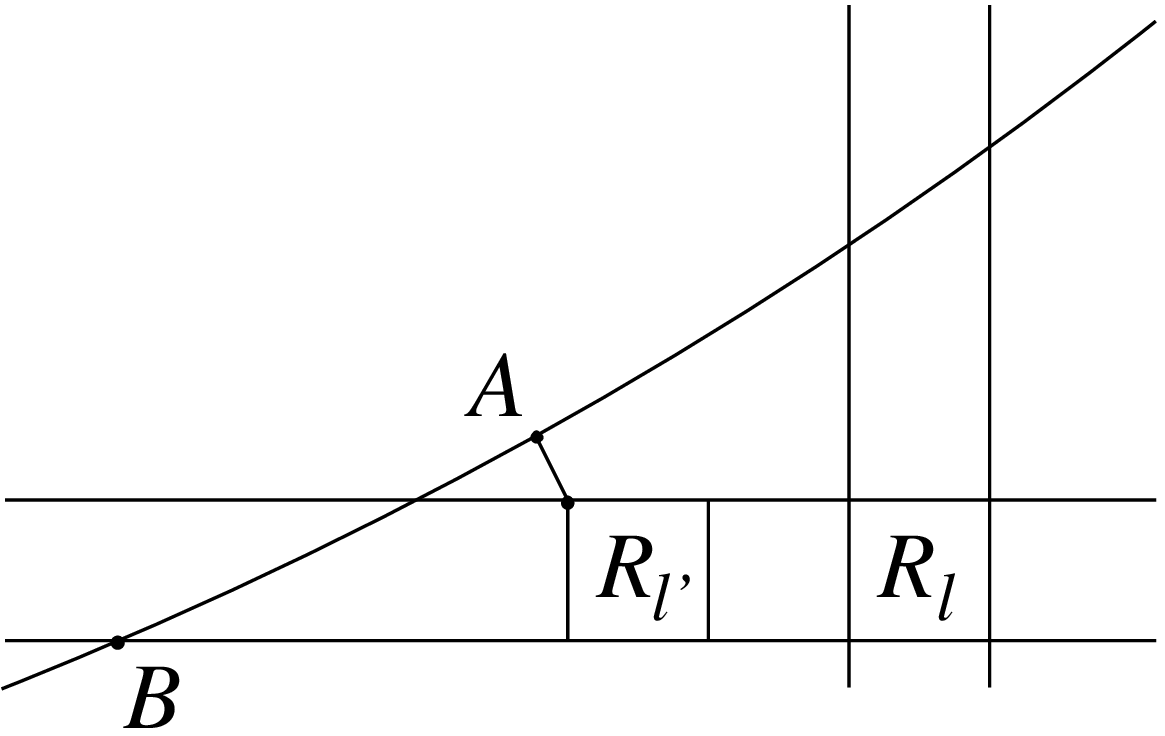}
}}
\\
\sc\figurename\ 5.1
\end{center}

Consider the case of intersecting $y$-projections (the other case is
similar). Then $R_{l'}$ is contained in
in the horizontal strip passing through $R_l$ (see Fig.~5.1).
By dilating along the $y$-axis, we may assume that $L=1$.
Since $\dist(R_{l'},W)\approx 2^{-m_l}$, there exists a point $A$ on
the graph of one of the functions $f_i(x)$ such that
$\dist(R_{l'},A)\approx 2^{-m_l}$. Let $B$ denote the point where the graph
of
$f_i(x)$ intersects the bottom of the strip. Since $f'_i(x)\approx L=1$,
we have $\dist(A,B)\lesssim  2^{-m_l}$, and therefore
$\dist(R_{l'},B)\lesssim 
 2^{-m_l}$. Thus all possible rectangles $R_{l'}$ are situated at a
distance $\lesssim 2^{-m_l}$ from no more than $n_\nu'$ points where the
bottom of the horizontal strip intersects $W$. 
This implies that the number of $R_{l'}$ in the horizontal strip is
bounded by a fixed constant, and the almost orthogonality is verified. 

Now let $R_l^*=(1+\eps)R_l$, where an $\eps>0$ is chosen so small that
$\dist(R_l^*,W)\approx 2^{-m_l}$ (in the anisotropic
norm). Consider a smooth partition of unity $\sum_l\gf_l=1$ on
$\bigcup R_l$ with $\supp \gf_l\subset R_l^*$, satisfying the natural
differential inequalities. We are going to decompose $T_{jk}$ using
this partition of unity. However, this decomposition will not be 
useful near the real multiple branches of $\cal Z$, since we will not have good
control on the size of $F$ there. For now we are just going to
localize away from those branches in the following way.

Let $\gb _i$ denote the power exponent of the first nonzero term in the 
asymptotic expansion of $\Im Y_{\nu i}(x)$; $\gb _i:=\infty$ if this
expansion is zero. For a large fixed number $Q$ we introduce the set 
\begin{align*}
W_Q=\mathop{{\bigcup}^*}_{i=1}^{n_\nu'}
\bigl\{(x,y)\in R_{jk}\bigl|
|y-f_i(x)|\<2^{-jQ}\bigr\},
\end{align*} 
where * indicates that the union is taken over all $i$ such that
$\gb_i=\infty$ and $f_i(x)\ne Y_{\nu i}(x)$. By the choice of
$f_i(x)$, this may happen only if the series $\bar Y_{\nu i}(x)$ is
real and there are several $\bar Y_{\nu i'}(x)$ having $\bar Y_{\nu i}(x)$ as
their asymptotic expansion. One can say that $W_Q$ is the tubular 
neighborhood of width $2^{-jQ}$ of the real multiple branches of $\cal
Z$ (see Fig.~1.4). 

The purpose of introducing $W_Q$ is that  
on $R_{jk}\backslash W_Q$ we have (if $j$ is large 
enough, which can be achieved by a further contraction of $V$)
\begin{equation}
\label{yminusfi}
|y-Y_{\nu i}(x)|\approx |y-f_i(x)|+2^{-j\gb_i}. 
\end{equation}  
Now let $\chi_Q$ be a smooth cut-off supported in the double of $W_Q$,
$\chi_Q\equiv 1$ on $W_Q$. We consider the decomposition
\begin{align}
\label{TQ}
&T_{jk}=T_Q+T^Q,\\ 
&T_Qf(x)=\int e^{i\gl
S(x,y)}\chi_Q(x,y)\chi_j(x)\chi_k(y)\chi(x,y)f(y)\,dy,\notag\\
&T^Q=\sum_l T^Q_l,\notag\\
&T^Q_lf(x)=\int e^{i\gl
S(x,y)}\gf_l(x,y)(1-\chi_Q(x,y))\chi_j(x)\chi_k(y)\chi(x,y)f(y)\,dy,\notag
\end{align}
In the rest of this section we prove that
$\|T^Q\|\lesssim\gl^{-\gd/2}$. The operator $T_Q$ will be dealt with
in the next section.

Let $T_l^Q$ be one of the operators from the decomposition of $T^Q$, and
assume that $T_l^Q\ne 0$, i.e.~that $R_l^*\cap
(R_{jk}\backslash W_Q)\ne\emptyset$.
Fix a point $(x_l,y_l)$ in this last intersection. We claim that 
\begin{equation}
\label{xlyl}
|y-f_i(x)|\approx |y_l-f_i(x_l)|
\end{equation}
for $(x,y)\in R_l^*$ and $i=1,\ldots,n_\nu'$. Indeed, let $(x',y')$
and $(x'',y'')$ be points of $R_l^*$ for which the value of
$|y-f_i(x)|$ is respectively minimal and maximal. Then
\begin{align*}
|y''-f_i(x'')|&\<|y'-f_i(x')|+|y''-y'|+|f_i(x'')-f(x')|\\
&\lesssim |y'-f_i(x')|+L2^{-m_l}\\
&\lesssim |y'-f_i(x')|,
\end{align*}
since $L^{-1}|y'-f_i(x')|\> \dist(R_l^*,W)\gtrsim 2^{-m_l}$. From this
\rf{xlyl} follows.

Now from \rf{fonrjk} and \rf{yminusfi} we see that on $R_l^*$
\[
|F|\approx 2^{-j(\gg_\nu B_\nu+A_\nu-\gg_\nu n_\nu')}
\prod_{i=1}^{n_\nu'}(|y_l-f_i(x_l)|+2^{-j\gb_i}) =:\mu_l.
\]
It follows by Lemma \ref{operatorcorput} (the condition (4) needs to be
checked, but this is easy) that $\|T^Q_l\|\lesssim (\mu_l\gl)^{-1/2}$.

We can get a lower bound on $\mu_l$ by noting that
$|y_l-f_i(x_l)|\gtrsim L2^{-m_l}=2^{-j(\gg_\nu-1)-m_l}$.
This gives
\[
\mu_l\gtrsim 2^{-j(\gg_\nu B_\nu+A_\nu)} 2^{(m_l-j)n_\nu'},
\]
and therefore
\begin{equation}
\label{tqlosc}
\|T^Q_l\|\lesssim \gl^{-1/2} 2^{j(\gg_\nu B_\nu+A_\nu)/2} 2^{-(m_l-j)n_\nu'/2}.
\end{equation}
On the other hand, by Corollary \ref{size},
\begin{equation}
\label{tqlsize}
\|T^Q_l\|\lesssim  2^{-m_l-j(\gg_\nu-1)/2}.
\end{equation}
It is clear that the passage from $R_l$ to $R_l^*$ preserved the
almost orthogonality of rectangles with fixed $m_l$. Therefore
\rf{tqlosc} and \rf{tqlsize} imply 
\[
\biggl\|\sum_{m_{l'}=m_l}T_{l'}^Q\biggr\|\lesssim\min\Bigl( 2^{-m_l-j(\gg_\nu-1)/2},
\gl^{-1/2} 2^{j(\gg_\nu B_\nu+A_\nu)/2} 2^{-(m_l-j)n_\nu'/2}\Bigr).
\]
Simple size considerations show that only rectangles with $m_l\>j-C$
may occur in the decomposition of $T^Q$. Therefore,
\begin{align*}
\|T^Q\|&\<\sumi_{m_l=j-C}\biggl\|\sum_{m_{l'}=m_l}T_{l'}^Q\biggr\|\\
&\lesssim \min\Bigl( 2^{-j(\gg_\nu+1)/2},
\gl^{-1/2} 2^{j(\gg_\nu B_\nu+A_\nu)/2}\Bigr).
\end{align*}
This is a familiar expression, and by taking the geometrical mean
killing the $j$-factor we obtain $\|T^Q\|\lesssim 
\gl^{-\gd_\nu/2}\<\gl^{-\gd/2}$.


\section{Estimates near the multiple real branches of $\cal Z$.}

To finish the proof of the theorem, we must estimate the operator
$T_Q$ appearing in the decomposition \rf{TQ} of $T_{jk}$. 

In the estimates below we can assume that $\gg_\nu\>1$, since this can
be achieved by passing to the adjoint operator if necessary.

Further, we can assume that $Q$ is chosen so large that the branches
of $\cal Z$ having different asymptotic expansions become completely separated in
the definition of $W_Q$. Since such branches can be treated separately,
we are reduced to the case when $T_Q$ has the form
\begin{eqnarray*}
T_Qf(x)=\int e^{i\gl S(x,y)}\chi_{jkQ}(x,y)f(y)\,dy,\\
\chi_{jkQ}(x,y)=\chi_j(x)\chi_k(y)\gf(2^{jQ}(y-g(x))).
\end{eqnarray*}
Here $\gf(t)$ is a $\cin$ cut-off supported in $[-1,1]$,
$g(x^{n!})\in\rr\cin$, $g(x)=cx^{\gg_\nu}+\ldots$, $c\ne0$, and in the
factorization \rf{factor1} exactly $N\>2$ functions $Y_{\nu i}(x)$
have asymptotic expansion coinciding with that of $g(x)$. We will
assume that this happens for $i=1,\ldots,N$. We also re-denote
$W_Q=\{(x,y)\in R_{jk}||y-g(x)|\<2^{-jQ}\}$.

We write $F(x,y)$ as
\[
F(x,y)=\tilde U(x,y)P(x,y),
\]
where $P(x,y)=\prod_{i=1}^N(y-Y_{\nu i}(x))$, and $\tilde U(x,y)$ is
the product of the rest of the terms in \rf{factor1}. 

Since all the branches of $\cal Z$ appearing in $\tilde U(x,y)$ 
are well separated from $W_Q$, there exists a constant $M_1\>0$ such
that
\[
|\tilde U|\approx 2^{-jM_1}\quad\text{on}\quad W_Q.
\]
Moreover, it can be seen directly that if $\sxy$ is completely
degenerate, we have $M_1=0$.

Further, by Proposition \ref{factor} (4), (5) we know that
$P(x,y)\in\rr\cin_+\(x^{1/n!}\)[y]$, so that $P(x,y)$ is $\cin$ in both
variables on $W_Q$. It is clear that 
\begin{equation}
\label{dyn}
\D_y^N P(x,y)=const\ne0. 
\end{equation}
We claim that, more generally,
\begin{equation}
\label{dxdy}
\D_x^k\D_y^{N-k} P(x,y)\ne 0\quad\text{on}\quad W_Q,\quad
k=0,\ldots,N.
\end{equation}
Denote $Q(x,y)=P(x^{n!},y)\in\cin_+\(x\)[y]$. The Taylor series of
$Q(x,y)$ is
\[
\bar Q(x,y)=\prod_{i=1}^N(y-\bar G(x)),\qquad \bar G(x)=\bar
g(x^{n!}).
\]
It is clear that
\begin{align*}
&[\D_x^l\D_y^{N-k}\bar Q](x,\bar G(x))=0,\qquad0\<l<k,\\
&[\D_x^k\D_y^{N-k}\bar Q](x,\bar G(x))\ne0.
\end{align*}
Therefore the factorizations of $\D_x^l\D_y^{N-k} Q(x,y)$, $l<k$, which can
be obtained as described in the proof of Lemma \ref{pfactor}, will
contain branches with the asymptotic expansion $\bar G(x)$, while the
factorization of $\D_x^k\D_y^{N-k} Q(x,y)$ will not contain such
branches. This implies \rf{dxdy}, provided that $Q$ is large enough,
since $\D_x^k\D_y^{N-k} P(x,y)$ can be expressed as 
\[
(\D_x^k\D_y^{N-k} Q)(x^{1/n!},y)+\sum_{l<k}c_l(x) (\D_x^l\D_y^{N-k}Q)(x^{1/n!},y)
\]
with coefficients $c_l(x)$ growing power-like as $x\to0$.

In addition, the above argument gives an estimate
\begin{equation}
\label{dxn}
\D_x^N P(x,y)\>2^{-jM_2}\quad\text{on}\quad W_Q,
\end{equation}
for some constant $M_2\>0$; $M_2=0$ if $\sxy$ is completely
degenerate.

Denote $\gs_j(x,y)=\frac 1{j!}\D^j_y P(x,y)$. Consider the
decomposition
\begin{align*}
&T_Q=\sumi_{l=-C}T_l,\\
&T_lf(x)=\int e^{i\gl
S(x,y)}\chi_{jkQ}(x,y)\bar\chi_l(\gs_0(x,y))f(y)\,dy,
\end{align*}
where $\bar\chi_l(t)$ is the characteristic function of the set
$2^{-l}\<|t|\<2^{-l+1}$, $C$ is a constant.

We are going to prove the estimates:
\begin{align}
\label{tl1}
&\|T_l\|\lesssim 2^{-{l}/{N}+{jM_2}/{2N}},\\
\label{tl2}
&\|T_l\|\lesssim \gl^{-1/2}(\log\gl)^{1/2}2^{l/2}l^{N-1/2}2^{jM_1/2}.
\end{align}
The required bound for $T_Q$ can then be derived as follows. 

Consider first the completely degenerate case, when 
$M_1=M_2=0$. We have
\[
\|T_Q\|\lesssim\sumi_{l=0}\min\bigl(2^{-l/N},\gl^{-1/2}2^{l/2}(\log\gl)^{1/2}l^{N-1/2}\bigr).
\]
If it were not for the factor of $(\log\gl)^{1/2} l^{N-1/2}$, the terms in
parentheses would become equal for $l=l_0=\frac N{N+2}\log \gl$, and we
would have the estimate $\|T_Q\|\lesssim \gl^{-1/(N+2)}$. In the
present situation we put $l_0=\frac N{N+2}\log \gl-k\log\log\gl$ with
indeterminate $k$ and have the estimate
\begin{align*}
\|T_Q\|&\lesssim 2^{-l_0/M}+\gl^{-1/2}2^{l_0/2}(\log\gl)^{1/2}l_0^{N-1/2}\\
&\lesssim
\gl^{-1/(N+2)}\Bigl[(\log\gl)^{k/N}+(\log\gl)^{N-k/2}\Bigr].
\end{align*}
The optimal value of $k$ is $k=\frac{2N^2}{N+2}$, which gives
\[
\|T_Q\|\lesssim \gl^{-\frac1{N+2}}(\log\gl)^{\frac{2N}{N+2}}
\]
in complete accordance with what is claimed in the theorem.

Assume now that $\sxy$ is not completely degenerate. In this case the
above argument gives in any case the estimate
\[
\|T_Q\|\<C_\eps 2^{jM}\gl^{-\frac1{N+2}+\eps}
\] 
for any $\eps>0$, with some constant $M$. (We do not pursue the possibility of
obtaining a $\log$ factor here, since as we will see in a moment, what we
have is already good enough.) Further, by Corollary \ref{size} we
certainly have the estimate
\[
\|T_Q\|\lesssim2^{-jQ/2}.
\]
The idea is that now we can take the geometric mean of the last two
estimates killing the $j$-factor and, if $Q$ is very large, this will
introduce only a very small increase in the exponent of $\gl$,
actually tending to zero as $Q\to\infty$. Thus we have 
\[
\|T_{Q_\eps}\|\<C_\eps\gl^{-\frac1{N+2}+\eps}.
\]  
We will show in a moment that in the case under consideration
$1/(N+2)>\gd/2$. 
This allows us to choose and fix $Q$ from the very beginning so large that 
$\|T_{Q}\|\lesssim\gl^{-\gd/2}$, thus proving the theorem.

We show that in fact $1/(N+2)>\gd_\nu/2$. Indeed, since we already
have $N$ branches whose expansion starts with $cx^{\gg_\nu}$, we know
that $n_\nu\>N$. Therefore $A_\nu=n_\nu\gg_\nu+A_\nu'\>N\gg_\nu$, and 
\[
\delta_\nu\<\frac{1+\gg_\nu}{1+N\gg_\nu+\gg_\nu}\<\frac 2{N+2},
\]
since $\gg_\nu\>1$. 
Besides that, the equality holds if and only if $\gg_\nu=1$,
$A_\nu=N$, $B_\nu=0$. But this corresponds exactly to the completely
degenerate case, which is excluded.

We now turn to the proof of the claimed bounds for $T_l$. The proof of
\rf{tl1} is easy and is based on the following well-known 

\begin{Lem}\tup{Christ \cite{Chr}, Lemma 3.3}
Let $f\in C^N[a,b]$ be such that $f^{(N)}\>\mu>0$ on $[a,b]$. Then for
any $\gg>0$
\[
\bigl|\{x\in [a,b]:|f(x)|\<\gg\}\bigr|\<A_k(\gg\mu)^{1/k},
\]
where the constant $A_k$ depends only on $k$.
\end{Lem}

By this lemma, in view of \rf{dyn} and \rf{dxn}, the kernel of $T_l$
is supported in a $y$-set of measure $\lesssim 2^{-l/N}$ for each
$x$, and in an $x$-set of measure $\lesssim 2^{-l/N+jM_2/N}$ for each
$y$. Now \rf{tl1} follows by Corollary \ref{size}.

The proof of \rf{tl2} is more intricate and is carried out by
a variation of a method developed in Seeger \cite{See1}, Section 3.
The key idea is to take an additional
dyadic localization in $\gs_j$, $1\<j\<N-1$. Let $l$ be fixed;
all constants below will however be independent of $l$.  Let
$\gg=(\gg_1,\ldots,\gg_{N-1})$ be a vector with integer components
$-C\<\gg_i\<l$, $C$ some constant. Denote
\[
\chi_\gg(x,y)=\chi_{jkQ}(x,y)\bar\chi_l(\gs_0(x,y))\prod_{i=1}^{N-1}
\dbarchi_{\gg_i}(\gs_i(x,y)),
\]
where $\dbarchi_{\gg_i}(t)$ is the characteristic function of the
set $2^{-\gg_i}\<|t|\<2^{-\gg_i+1}$ for $\gg_i<l$, and of the
set $|t|\<2^{-l+1}]$ for $\gg_i=l$.

For an appropriate fixed $C$ we have the decomposition
\begin{align*}
&T_l=\sum_\gg T_\gg,\\
&T_\gg f(x)=\int e^{i\gl S(x,y)}\chi_\gg(x,y)f(y)\,dy.
\end{align*}
We are going to prove that for each $\gg$
\begin{equation}
\label{tgamma}
\|T_\gg\|\lesssim \gl^{-1/2}(\log\gl)^{1/2}2^{l/2}l^{1/2}2^{jM_1/2}.
\end{equation}
This will imply \rf{tl2}, since the number of $T_\gg$ in the
decomposition of $T_l$ is $\lesssim l^{N-1}$.
 
The kernel of the operator $T_\gg^*T_\gg$ has the form
\[
K(y_2,y_1)=\int
e^{i\gl[S(x,y_2)-S(x,y_1)]}\chi_\gg(x,y_1)\chi_\gg(x,y_2)\,dx.
\]
Assuming that $y_2>y_1$, and using Taylor's formula in $y$ for
$P(x,y)$, we have 
\begin{align}
\label{[]}
[S(x,y_2)-S(x,y_1)]'_x&=\int_{y_1}^{y_2}\tilde
U(x,y)P(x,y)\,dy\\
\notag&=\int_{y_1}^{y_2}\tilde U(x,y)\Bigl[\sum_{j=0}^N
\gs_j(x,y_1)(y-y_1)^j\Bigr]\,dy\\
\notag&=\sum_{j=0}^N\gs_j(x,y_1)\int_{y_1}^{y_2}\tilde U(x,y)(y-y_1)^j\,dy.
\end{align}
Notice that $\int_{y_1}^{y_2}\tilde U(x,y)(y-y_1)^j\,dy\approx
2^{-jM_1}(y_2-y_1)^{j+1}$. So the RHS of \rf{[]} looks like a
polynomial in $y_2-y_1$ with dyadically restricted coefficients. To
handle such polynomials, we need the following variant of Lemma 3.2
from \cite{See1}. We chose to give a proof, since we have found a
one simpler than in \cite{See1}.

\begin{Lem} 
\label{dyadpol}
For an integer $N\>1$, an integer vector
$r=(r_1,\ldots,r_N)$, $r_i\>0$, and a constant $C>0$ consider the set
$\cal P=\cal P(r,C,N)$ of all polynomials of the form
$
P(h)=1+\sum_{i=1}^Na_ih^i
$
with real coefficients $a_i$ satisfying 
\begin{align*}
&|a_i|\in[C^{-1}2^{r_i},C2^{r_i}]\quad\text{if}\quad r_i>0,\\
&|a_i|\<C\qquad\qquad\quad\quad\text{if}\quad r_i=0.
\end{align*}
Then there exists a constant $B=B(C,N)$, independent of $r$, and a set
$E\in[0,1]$ of the form 
\begin{equation}
\label{E}
E=[0,2^{\gb_1}]\cup[2^{\ga_2},2^{\gb_2}]\cup\ldots\cup[2^{\ga_s},2^{\gb_s}],
\end{equation}
such that\\ 
\tup{1} $\ga_i$, $\gb_i$ are integers, $\gb_1<\ga_2<\gb_2<\ldots<\ga_s<\gb_s\<0$,\\
\tup{2} $s\<B${\rm ;} $\gb_1\>-B\max(r_i)${\rm;}
$\bigl[(1-\gb_s)+\sum_{j=1}^{s-1}(\ga_{j+1}-\gb_j)\bigr]\<B$,\\
\tup{3} $|P(h)|\>B^{-1}$ for $h\in E$ for any $P\in\cal P$.
\end{Lem}
\begin{proof} Put $r_0=0$.
Consider the convex set $\GS$ given as the intersection of the
half-planes lying above the lines $y=r_i+ix$,
$i=0,\ldots,N$. Let $A_i$, $i=1,\ldots,n$ be all the corner points of
the boundary of $\GS$ with the $x$-coordinates
$x_1<x_2<\ldots<x_n$. It is clear that $n\<N$. 
We claim that for any $P\in \cal P$
\[
|P(h)|\>B^{-1}\quad\text{if}\quad h>0,\quad\log
 h\notin\bigcup_{j=1}^n(x_j-B,x_j+B).
\]
It is not difficult to see that this implies (1)--(3).

Let $\log h\in[x_j+B,x_{j+1}-B]$, and assume that the boundary points
$A_j$ and $A_{j+1}$ belong to the line $y=r_k+kx$. Since
$A_j$ and $A_{j+1}$ lie above all the other lines $y=r_i+ix$, we have for
all $i$
\begin{align*}
r_i+ix_j&\<r_k+kx_j,\\
r_i+ix_{j+1}&\<r_k+kx_{j+1}.
\end{align*}
From these two estimates it follows that
\begin{align*}
|a_ih^i|&\lesssim |a_kh^k|2^{(k-i)(\log h-x_j)},\\
|a_ih^i|&\lesssim |a_kh^k|2^{(i-k)(x_{j+1}-\log h)}.
\end{align*}
All in all,
\[
|a_ih^i|\lesssim |a_kh^k|2^{-|k-i|B}.
\]
This estimate clearly implies $|P(h)|\gtrsim|a^kh^k|\gtrsim2^{kB}\>1$,
provided that $B$ is large enough.
\end{proof}

Now if we take out the factor of
$2^{-l-jM_1}(y_2-y_1)$, the expression in the RHS of \rf{[]} has the form of polynomial in $h=y_2-y_1$
falling under the scope of the lemma with $r_i=l-\gg_i$. So we have a
set $E$ of the form \rf{E} such that
\[
\bigl|[S(x,y_2)-S(x,y_1)]'_x\bigr|\gtrsim 2^{-l-jM_1}(y_2-y_1)\quad
 \text{if}\quad y_2-y_1\in E.
\]
We claim that this implies
\begin{equation}
\label{|K|}
|K(y_2,y_1)|\lesssim2^{l+jM_1}\gl^{-1}(y_2-y_1)^{-1}\qquad(y_2-y_1\in
 E).
\end{equation}
Indeed, this will follow from Lemma \ref{vandercorput} with $k=1$, if
we prove that there exists a constant $C$ independent of $y_1$ and
$y_2$ such that for fixed $y_1$ and $y_2$\\
(1) the number of intervals of monotonicity of
$[S(x,y_2)-S(x,y_1)]'_x$ considered as a function of $x$ is less than
$C$,\\
(2) the number of intervals comprising the $x$-set where
$\chi_\gg(x,y_1)\chi_\gg(x,y_2)$ is non-zero is less than $C$.

  To show (1), note that $\D_x^N F(x,y)\ne 0$ on $W_Q$. It follows that
\[
\D_x^{N+1}[S(x,y_2)-S(x,y_1)]=\int_{y_1}^{y_2}\D_x^N F(x,y)\,dy\ne 0.
\]
Therefore, $[S(x,y_2)-S(x,y_1)]''_{xx}$ vanishes at most $N-1$ times,
which implies (1).

To show (2), it suffices to check that the number of intervals in the
set $\{x|(x,y)\in W_Q,\ a\<\gs_j(x,y)\<b\}$ is bounded by a constant
independent of $a$ and $b$ for each $0\<j\<N-1$. However, this last
statement follows from \rf{dxdy}.

Unfortunately, to prove the claimed norm estimate for $T_\gg$, we will
need still another decomposition taking into account the form of the
set $E$. Namely, for $1\<k\<s$ and an integer $n$ we put
\[
\chi_{kn}(y)=\gp(2^{\gb_k}y-n), 
\]
where $\gp(t)$ is the characteristic function of the interval $[0,1]$,
and consider the operators 
\[
T_{kn}f(x)=\int e^{i\gl S(x,y)}\chi_\gg(x,y)\chi_{kn}(y)f(y)\,dy.
\]
We are going to prove by induction in $k$ that for each $n$
\[
\|T_{kn}\|\lesssim \gl^{-1/2}(\log\gl)2^{l/2}l^{1/2}2^{jM_1/2}.
\]
The statement for $k=s$ implies the required estimate \rf{tgamma},
since $T_\gg=\sum_n T_{sn}$, and the sum contains no more than
$2^{-\gb_s}\<C$ terms.

For $k=1$, we use the kernel of the operator $T_{1n}^*T_{1n}$, which has the form
\[
\chi_{1n}(y_1)\chi_{1n}(y_2)K(y_2,y_1),
\]
where $K(y_2,y_1)$ is the kernel of $T_\gg^*T_\gg$. If this expression
is not zero, then $|y_2-y_1|\<2^{\gb_s}$. In view of \rf{|K|}, and
also because $|K|\lesssim 1$, Lemma \ref{schur} says that
\begin{align*}
\|T_{1n}^*T_{1n}\|&\lesssim
\int_0^{2^{\gb_s}}\min(1,2^{l+jM_1}\gl^{-1}t^{-1})\,dt\\
&\lesssim2^{l+jM_1}\gl^{-1}\int_0^{\gl2^{\gb_s-l-jM_1}}\min(1,t^{-1})\,dt\\
&\lesssim2^{l+jM_1}\gl^{-1}\log\gl,
\end{align*}
which is even better by a factor of $l$ than what we need.

The induction step is performed by using the decomposition
\[
T_{k+1,n}=\sum_{n'}T_{kn'}.
\]
We will need the following variant of the Cotlar--Stein lemma, which
can be proved by an easy adaptation of the standard proof given in
\cite{STEIN}, see e.g.~Comech \cite{Comech}, Appendix.

\begin{Lem} 
\label{cotlarstein}
Let $T_i$ be a family of operators on a Hilbert space
$H$ such that \\
\tup{1} $T_iT_{i'}^*=0$ for $i\ne i'$,\\
\tup{2} $\sum_{i'}\|T_i^*T_{i'}\|\<C$ with a constant $C$ independent of
$i$.\\
Then $\|\sum T_i\|\<C^{1/2}$.
\end{Lem}  

We have $T_{kn'}T_{kn''}^*=0$ for $n'\ne n''$. Let us estimate the sum 
\begin{equation}
\label{sum}
\sum_{n''}\|T_{kn'}^*T_{kn''}\|
\end{equation} 
for a fixed $n'$.  
Since both $T_{kn'}$ and $T_{kn''}$ appear in the decomposition of $T_{k+1,n}$,
we have $|n'-n''|\<2^{\gb_{k+1}-\gb_k}$. Further, the kernel of
$T_{kn'}^*T_{kn''}$ has the form
\[
\chi_{kn'}(y_1)\chi_{kn''}(y_2)K(y_2,y_1).
\]
If this expression is different from zero, then
\begin{equation}
\label{yy}
2^{\gb_k}|y_2-y_1|\in[|n'-n''|-1,|n'-n''|+1].
\end{equation}
Assume first that 
\begin{equation}
\label{nn}
2^{\ga_{k+1}-\gb_{k}}+1\<|n'-n''|\<2^{\gb_{k+1}-\gb_{k}}-1.
\end{equation} 
Then \rf{yy} implies $|y_2-y_1|\in E$, and we can use the estimate
\rf{|K|}. By Lemma \ref{schur},
\begin{align*}
\|T_{kn'}^*T_{kn''}\|&\lesssim\int_{2^{-\gb_k}(|n'-n''|-1)}^{2^{-\gb_k}(|n'-n''|+1)}
2^{l+jM_1}\gl^{-1}t^{-1}\,dt\\
&\lesssim 2^{l+jM_1}\gl^{-1}|n'-n''|^{-1}.
\end{align*}
Therefore the part of the sum \rf{sum} over $n''$ satisfying \rf{nn}
is bounded by
\[
2^{l+jM_1}\gl^{-1}\sum_{m=2^{\ga_{k+1}-\gb_{k}}}
^{2^{\gb_{k+1}-\gb_{k}}}\frac 1m\lesssim 2^{l+jM_1}\gl^{-1} 
(\gb_{k+1}-\ga_{k+1})\<2^{l+jM_1}\gl^{-1}l,
\]
where we used the fact that by Lemma \ref{dyadpol} (2) $\gb_1\>-B l$.

However, the number of $n''$ which do not satisfy \rf{nn} is bounded
by a constant in view of Lemma \ref{dyadpol} (2), so the
corresponding part of \rf{sum} is bounded by
$C\sup\|T_{kn''}\|^2\lesssim l2^{l+jM_1}(\log\gl)\gl^{-1}$ by the induction 
hypothesis. 

By applying Lemma \ref{cotlarstein}, we complete the induction
step. The theorem is now proved.

\bibliographystyle{amsplain}

\end{document}